\pdfoutput=1

\documentclass[11pt]{article}

\usepackage{graphicx,color,geometry}
\usepackage{amssymb,amsfonts,amsmath,amsthm,cite}
\newcommand{\eps}{\varepsilon}
\newtheorem{theorem}   {Theorem}
\newtheorem{lemma}   {Lemma}
\newtheorem{remark}   {Remark}

\begin{document}

\title{\bf Self-similarity and long-time behavior of solutions of the
  diffusion equation with nonlinear absorption and a boundary source}

\author{ Peter V. Gordon\thanks{Department of Mathematical Sciences,
    New Jersey Institute of Technology, Newark, NJ 07102, USA} \and
  Cyrill B. Muratov$^*$}
    \maketitle

    {\em Dedicated to Hiroshi Matano on the occasion of his 60th
      birthday.}

\begin{abstract}
  This paper deals with the long-time behavior of solutions of
  nonlinear reaction-diffusion equations describing formation of
  morphogen gradients, the concentration fields of molecules acting as
  spatial regulators of cell differentiation in developing
  tissues. For the considered class of models, we establish existence
  of a new type of ultra-singular self-similar solutions. These
  solutions arise as limits of the solutions of the initial value
  problem with zero initial data and infinitely strong source at the
  boundary. We prove existence and uniqueness of such solutions in the
  suitable weighted energy spaces. Moreover, we prove that the
  obtained self-similar solutions are the long-time limits of the
  solutions of the initial value problem with zero initial data and a
  time-independent boundary source.
\end{abstract}

\section{Introduction}

In the studies of reaction-diffusion equations, one canonical problem
deals with the following equation \cite{galaktionov85,brezis83}:
\begin{align}
  \label{eq:pde}
  u_t = \Delta u - u^p, \qquad (x, t) \in \mathbb R^d \times (0,
  \infty).
\end{align}
Here $p > 1$ is a constant and $u = u(x, t) > 0$ can be viewed as the
concentration of a chemical species diffusing in the $d$-dimensional
space subject to degradation whose rate is an increasing function of
the species concentration. Usually, one considers the associated
Cauchy problem with some non-negative initial data $u(x, 0) =
u_0(x)$. During the 1980's, this problem attracted a considerable
attention, in particular in the case of measure-valued initial data
(e.g., when $u_0$ is a Dirac mass)
\cite{brezis83,gmira84,brezis86,kamin85,escobedo87,oswald88}. In the
course of these studies, it was discovered that \eqref{eq:pde} possess
self-similar solutions for all $1 < p < (2 + d) / d$, which are smooth
for all $t > 0$ and converge to zero outside the origin, while blowing
up at the origin when $t \to 0^+$ \cite{galaktionov85,brezis86} (see
also \cite{escobedo87} for a variational approach). These solutions
play important roles in determining the long-time behavior of the
solutions of the Cauchy problem for general classes of initial data
and in some sense describe the transient dynamics in systems described
by \eqref{eq:pde} \cite{galaktionov85, kamin85, oswald88, escobedo95,
  bricmont96, wayne97, herraiz99}. In particular, a special class of
self-similar solutions of \eqref{eq:pde} called {\em very singular
  solutions} attract the physically important class of initial data
with sufficiently fast asymptotic decay
\cite{kamin85,escobedo88,galaktionov85}.

Equation \eqref{eq:pde} with $p \geq 1$ on domains with boundaries
also arises as a canonical model of morphogen gradient formation (for
recent reviews, see \cite{rmss:dc06,llnw09,othmer09,wkg09}). Morphogen
gradients are concentration fields of molecules acting as spatial
regulators of cell differentiation in developing tissues
\cite{ms02}. In particular, the case $p > 1$ was proposed to describe
a robust patterning mechanism whereby morphogen increases the
production of molecules which, in turn, increase the rate of morphogen
degradation \cite{eldar03}. For example, a protein called Sonic
hedgehog (Shh) is known to induce the expression of its receptor
Patched, which both transduces the Shh signal and mediates Shh
degradation by cells in the {\em Drosophila} embryo
\cite{cs96,ilrekr00}. 

An important aspect of morphogen dynamics is the presence of localized
sources at the boundary of the morphogenetic field. This leads to the
need to consider initial boundary value problems, whose prototype is
the following one-dimensional problem:
\begin{eqnarray}
  \label{eq:u}
  \left\{\begin{array}{ll}
      u_t=u_{xx}-u^p  &  (x,t)\in[0,\infty)\times(0,\infty), \\
      u_x(0, t)=-\alpha &  t\in(0,\infty),\\
      u(x, 0)=0 & x\in [0,\infty).
\end{array}\right.
\end{eqnarray}
This problem can be viewed as an extension of the Cauchy problem for
\eqref{eq:pde} defined for $x > 0$ in the presence of a boundary
source at $x = 0$. Here $\alpha > 0$ is a constant characterizing the
source strength of morphogen production, and the zero initial
condition corresponds to the absence of the morphogen at the onset of
patterning. In what follows, we will restrict our attention only to
this simplest model of morphogen gradient formation.

In the context of morphogenesis, one is often interested in the
establishment of a stationary morphogen profile and the transient
dynamics that leads to it. The stationary problem for \eqref{eq:u} can
be written as the following boundary value problem:
\begin{eqnarray}\label{eq:v}
  v_{xx}-v^p=0, \qquad 
  v_x(0)=-\alpha, \qquad  v(\infty) = 0,
\end{eqnarray}
whose unique solution for any $p > 1$ is explicitly given by
\begin{eqnarray}
  \label{eq:a}
  v(x)= \left( \frac{2 (p+1)}{(p-1)^2} \right)^{1 \over p-1}
  (a+x)^{-\frac{2}{p- 1}}, \quad
  a = \left( \frac{2^{\frac{p}{p+1}} (p+1)^{\frac{1}{p+1}} }{p-1} \right)
  \alpha^{-\frac{p-1}{p+1}}.  
\end{eqnarray}
In fact, it is easy to see that the stationary solution $v(x)$ in
\eqref{eq:v} is the limit of the solution $u(x, t)$ of \eqref{eq:u} as
$t \to \infty$ for each $x \geq 0$, and is approached monotonically
from below \cite{gsbms:pnas11}. However, as we noted in
\cite{gsbms:pnas11}, this approach is not uniform in $x$ and for each
fixed $x \geq 0$ occurs on the diffusive time scale $\tau_p(x) =
O(x^2)$, which diverges as $x \to \infty$. Thus, the timing of the
establishment of the steady morphogen concentration at a given point
depends rather sensitively on the location of that point.

\begin{figure}[t]
  \centering
  \includegraphics[width=4in]{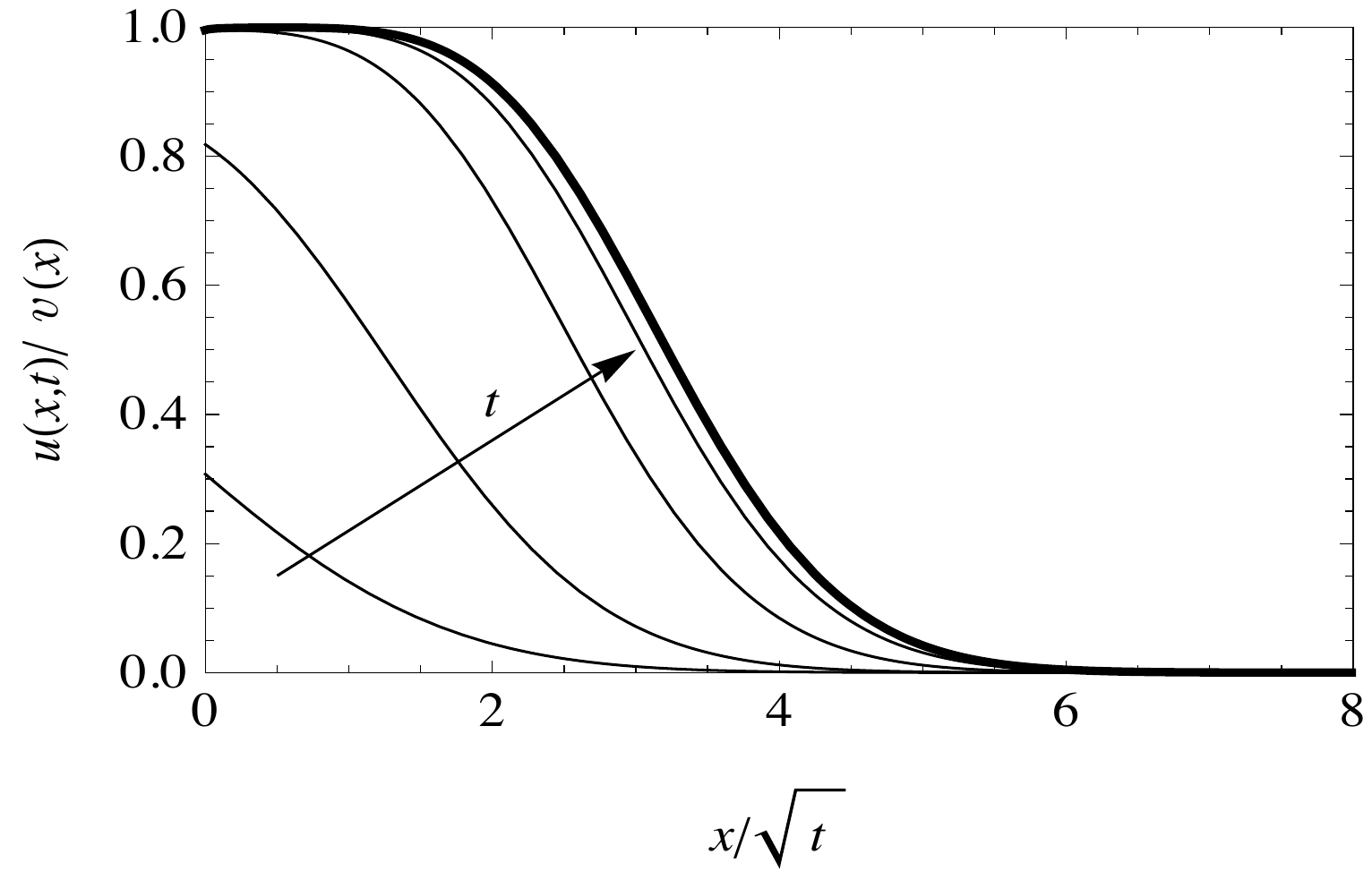}
  \caption{Numerical solution of \eqref{eq:u} in self-similar
    variables for $p = 2$ and $\alpha = 1$. Thin lines show snapshots
    of the solution corresponding to $t = 0.1, 1, 10, 100$ (the
    direction of time increase is indicated by the arrow). The bold
    line shows the asymptotic solution. \label{fig:collapse}}
\end{figure}

To better understand the dynamics of the approach of the solution of
\eqref{eq:u} to the stationary solution, we undertook numerical
studies of the initial boundary problem in \eqref{eq:u} for various
values of $p > 1$. In those studies, we discovered that when the ratio
of the solution at a given $x$ to the value of the stationary solution
at $x$ is plotted vs. the diffusion similarity variable $x /
\sqrt{t}$, the numerical solution approaches some universal limit
curve depending only on the value of $p$ \cite{mgs:pre11}. This
process is illustrated in Fig. \ref{fig:collapse}, where the results
are presented for the biophysically important case $p = 2$. This
observation suggested to us some hidden self-similarity in the
behavior of solutions of \eqref{eq:u} \cite{barenblatt}. A simple
scaling argument indicates that the long-time behavior of the solution
of \eqref{eq:u} for a fixed value of $\alpha > 0$ is closely related
to the behavior of solutions of \eqref{eq:u} at fixed $x > 0$ and $t >
0$ as $\alpha \to \infty$ \cite{mgs:pre11}. We found numerically that
in the limit $\alpha \to \infty$ the solutions of \eqref{eq:u} attain
a self-similar profile (see the following section for precise
definitions) \cite{mgs:pre11}. The purpose of this paper is to
substantiate these numerical observations by establishing existence
and properties of what we will call {\em ultra-singular self-similar
  solutions} in the limit of infinite boundary source strength. We
also prove that these solutions are indeed the long-time limits of the
solutions of \eqref{eq:u} in the above sense.

We note that the solutions constructed by us form a new class of
self-similar solutions to \eqref{eq:pde} in $d = 1$. Indeed, our
solutions can be trivially extended to the whole real line by a
reflection and can be viewed as singular solutions of \eqref{eq:pde}
that blow up at the origin. We point out that these solutions are
different from the self-similar solutions studied in
\cite{galaktionov85,brezis86}. The ultra-singular solutions of
\eqref{eq:pde} constructed by us can be viewed as the more singular
counterparts of the very singular solutions of \cite{brezis86} in the
following sense: the singularity in the former is concentrated on a
half-line ($x = 0, t > 0$) in the $(x, t)$ plane, while the
singularity in the latter occurs only at a single point ($x = 0, t =
0$). Similarly, our convergence result for the solutions of
\eqref{eq:u} with $\alpha \in (0, \infty)$ may be viewed as a
counterpart of the result of \cite{kamin85}, in the sense that in the
former case the solution can be viewed as the distributional solution
of \eqref{eq:pde} with an added term $2 \alpha \delta(x)$ in the
right-hand side, while in the latter case one can think of the
solution as the distributional solution of \eqref{eq:pde} with the
term $\alpha \delta(x) \delta(t)$ added to the right-hand side.

Before concluding this section, let us briefly mention a few possible
extensions and open problems related to our present work. It would be
interesting to understand the role our self-similar solutions play for
the singular solutions of the initial value problem associated with
\eqref{eq:pde} for general non-zero initial data. Let us point out
that even the basic questions of existence and uniqueness of such
singular solutions for the considered parabolic problems in suitable
function classes are currently open (see \cite{veron11} for a very
recent related work). Other natural extensions include higher
dimensional versions of the considered problem, as well as a proof of
global stability of self-similar solutions. These studies are
currently ongoing. From the point of view of applications, it is also
important to consider solutions of \eqref{eq:pde} with added
time-varying singular sources, for which both the very singular and
the ultra-singular solutions may be relevant.

Our paper is organized as follows. In
Sec. \ref{sec:similarity-ansatz}, we introduce a singular version of
the initial boundary value problem in \eqref{eq:u} and prove
existence, uniqueness, monotonicity and limiting behavior of the
self-similar solution to this singular problem. Then, in
Sec. \ref{s:long} we prove that the obtained self-similar solutions
are the long-time limits of the solutions of \eqref{eq:u} in an
appropriate sense.

\section{Singular solutions and the similarity ansatz}
\label{sec:similarity-ansatz}

Let us consider \eqref{eq:u} with {\em infinite source at the
  boundary}, i.e., the following singular initial boundary value
problem:
\begin{eqnarray}
  \label{eq:uu}
  \left\{\begin{array}{ll}
      u_t=u_{xx}-u^p  &  (x, t)\in(0,\infty)\times(0,\infty), \\
      u(0,t)=\infty,  &  t\in(0,\infty),\\
      u(x,0)=0 & x\in (0,\infty).
\end{array}\right.
\end{eqnarray}
By a solution to  \eqref{eq:uu}, we mean a classical solution for
all $(x, t) \in (0, \infty) \times (0, \infty)$ decaying sufficiently
fast as $x \to +\infty$ for all $t > 0$, and continuous up to $t = 0$
for all $x > 0$. Note that for each $p> 1$ this problem possesses a
{\em singular stationary solution}
\begin{eqnarray}
  \label{eq:vinf}
  v_\infty(x) = \left( \frac{2 (p+1) }{(p-1)^2} \right)^{1
    \over p - 1} \left( {1 \over x} \right)^{2 \over p - 1},
\end{eqnarray}
which is the limit of $v_\alpha(x)$ as $\alpha \to \infty$ for each $x
> 0$. 

Consistently with the discussion in the introduction, we now seek
solutions of \eqref{eq:uu} in the form
\begin{eqnarray}
  \label{eq:simanz}
  u(x, t) = v_\infty(x) \phi (\zeta), \quad \zeta=\ln(x / \sqrt{t}),
\end{eqnarray}
for some function $0 \leq \phi(\zeta) \leq 1$, which will be referred
to as the {\em self-similar profile}. Substituting the similarity
ansatz from \eqref{eq:simanz} into \eqref{eq:uu}, after some algebra
we obtain the following equation for the self-similar profile $\phi$:
\begin{eqnarray} 
  \label{eq:phiz} {d^2 \phi \over d \zeta^2} +
  \left(\frac{e^{2\zeta}}{2} - \frac{p+3}{p-1}\right){d \phi \over d
    \zeta} + 
  \frac{2(p+1)}{(p-1)^2}\phi(1-\phi^{p-1})=0,
\end{eqnarray}
which must hold for all $\zeta \in (-\infty, \infty)$, supplemented
with the limit behavior
\begin{eqnarray}
  \label{eq:phibcz1}
  \lim_{\zeta \to -\infty} \phi(\zeta) = 1, \qquad \lim_{\zeta \to
    -\infty} {d \phi(\zeta) \over d \zeta} = 0,  \\ 
  \lim_{\zeta \to +\infty} \phi(\zeta) = 0, \qquad
  \lim_{\zeta \to +\infty} 
  {d \phi(\zeta) \over d \zeta} = 0.   \label{eq:phibcz2}
\end{eqnarray}
Existence and multiplicity of solutions of \eqref{eq:phiz} satisfying
\eqref{eq:phibcz1} and \eqref{eq:phibcz2} are not at all {\em a
  priori} obvious in view of both the non-linearity and the presence
of singular terms in the considered boundary value problem. In
\cite{mgs:pre11}, we were able to construct such solutions numerically
for several values of $p$. Here we establish their existence and
uniqueness for all $p > 1$ within a natural class of functions.

We will prove existence and uniqueness of solutions of \eqref{eq:phiz}
satisfying \eqref{eq:phibcz1} and \eqref{eq:phibcz2} in the weighted
Sobolev space $H^1(\mathbb R, d \mu)$, which is obtained as the
completion of the family of smooth functions with compact support with
respect to the Sobolev norm $||.||_{H^1(\mathbb R, d \mu)}$, defined
as
\begin{eqnarray}
  \label{eq:H1}
  ||w||_{H^1(\mathbb R, d \mu)}^2  =  ||w_\zeta||_{L^2(\mathbb R, d
    \mu)}^2 +  ||w||_{L^2(\mathbb R, d \mu)}^2,
\end{eqnarray}
where $||w||_{L^2(\mathbb R, d \mu)}^2 = \int_\mathbb{R} w^2(\zeta)
d \mu(\zeta)$, and the measure $d \mu$ is
\begin{eqnarray}
  \label{eq:muu}
  d \mu(\zeta) = \rho(\zeta) d \zeta, \qquad \rho(\zeta) =  
  \exp\left\{\frac{e^{2\zeta}}{4}-\left(\frac{p+3}{p-1}\right) \zeta 
  \right\}.
\end{eqnarray}
Our existence and uniqueness result is given by the following
theorem. 

\begin{theorem}
  \label{t:exun}
  There exists a unique weak solution $\phi$ of  \eqref{eq:phiz},
  such that $\phi - \eta \in H^1(\mathbb R, d \mu)$, with $\mu$
  defined in  \eqref{eq:muu}, for every $\eta \in C^\infty(\mathbb
  R)$, such that $\eta(\zeta) = 1$ for all $\zeta \leq 0$ and
  $\eta(\zeta) = 0$ for all $\zeta \geq 1$. Furthermore, $\phi \in
  C^\infty(\mathbb R)$, satisfies  \eqref{eq:phiz} classically and
  $0 < \phi < 1$. In addition, $\phi$ is strictly decreasing and
  satisfies  \eqref{eq:phibcz1} and \eqref{eq:phibcz2}.
\end{theorem}

Before proceeding to the proof of Theorem \ref{t:exun}, let us
establish a basic technical lemma needed to deal with the weighted
spaces introduced above, which is an extension of \cite[Lemma
4.1]{lmn:cpam04} for exponentially weighted Sobolev spaces (cf. also
\cite[Lemma 1.5]{escobedo87}).

\begin{lemma}
  \label{l:poinc}
  Let $w \in H^1(\mathbb R, d \mu)$. Then there exists $R_0 > 0$ such
  that 
  \begin{eqnarray}
  \label{eq:poinc}
  \int_R^\infty w^2 d \mu \leq \frac{e^{-2 R}}{2} \int_R^\infty
  \left(\frac{d w}{d\zeta}\right)^2 d \mu \qquad \qquad \forall R \geq
  R_0, 
 \end{eqnarray}
and 
 \begin{eqnarray}
  \label{eq:poincR}
  \rho(R) w^2(R)\le 2 e^{-R} \int_R^{\infty}\left(\frac{d
      w}{d\zeta}\right)^2 d \mu \qquad \qquad \text{for a.e. } R \geq
  R_0. 
   \end{eqnarray}
Moreover, there exists $R_0' < 0$ such that
  \begin{eqnarray}
  \label{eq:poinc+}
  \int_{-\infty}^R w^2 d \mu \leq 8\left( \frac{p-1}{p+3}\right)^2 \int_{-\infty}^R
  \left(\frac{d w}{d\zeta}\right)^2 d \mu \qquad \qquad \forall R \leq
  R_0', 
 \end{eqnarray}
and 
\begin{eqnarray}
  \label{eq:poincR+}
  \rho(R) w^2(R)\le 8 \left( \frac{p-1}{p+3} \right)
  \int_{-\infty}^R\left(\frac{d  w}{d\zeta}\right)^2 d \mu \qquad
  \qquad \text{for a.e. } R \leq R_0'. 
\end{eqnarray}

\end{lemma}

{\em \bf \small \noindent Proof. } Arguing by approximation, observe
that by an explicit computation and an application of Cauchy-Schwarz
inequality we have
\begin{eqnarray}
  \frac12\left( w^2(R)\rho(R)+\int_R^{\infty} \left(\frac{d
      }{d\zeta}\ln \rho\right) w^2 d\mu \right) \qquad \qquad \qquad
  \qquad 
  \nonumber \\  
  =-\int_R^{\infty}
  w\frac{d w}{d \zeta}d\mu\le \left( \int_R^{\infty}
    w^2d\mu\int_R^{\infty}
    \left(\frac{dw}{d\zeta}\right)^2d\mu\right)^{1/2} .
\label{csrho}
\end{eqnarray}
In particular, \eqref{eq:muu} and \eqref{csrho} yield 
\begin{eqnarray}
\left(\frac{e^{2R}}{2}-\frac{p+3}{p-1}\right)^2
 \int_R^{\infty} w^2d\mu
\le4 \int_R^{\infty} \left(\frac{dw}{d\zeta}\right)^2d\mu, 
\end{eqnarray}
which for large enough $R$ implies \eqref{eq:poinc}.  Next, since ${d
  \over d \zeta} \ln \rho > 0$ for large positive $\zeta$, dropping
the second term in the left-hand side of \eqref{csrho} and using
\eqref{eq:poinc}, we obtain \eqref{eq:poincR}.
 
Similarly, we note that
\begin{eqnarray}
  \frac12\left( w^2(R)\rho(R)-\int_{-\infty}^R \left(\frac{d
      }{d\zeta}\ln \rho\right) w^2 d\mu \right) \qquad \qquad \qquad
  \qquad 
  \nonumber \\  
  =\int_{-\infty}^R w\frac{d w}{d \zeta}d\mu\le \left(
    \int_{-\infty}^R w^2d\mu\int_{-\infty}^R
    \left(\frac{dw}{d\zeta}\right)^2d\mu\right)^{1/2}, 
  \label{csrhoo}
\end{eqnarray}
which implies
\begin{eqnarray}
  \left(\frac{p+3}{p-1}-\frac{e^{2R}}{2}\right)^2
  \int_{-\infty}^R w^2d\mu
  \le4 \int_{-\infty}^R \left(\frac{dw}{d\zeta}\right)^2d\mu,
\end{eqnarray}
and thus \eqref{eq:poinc+} holds for sufficiently large negative $R$.
Finally since ${d \over d \zeta} \ln \rho < 0$ for large negative
$\zeta$, from \eqref{csrhoo} and \eqref{eq:poinc+} we obtain
\eqref{eq:poincR+}. \qed

\medskip

{\em \bf \small \noindent Proof of Theorem \ref{t:exun}. } The proof
consists of five steps. \smallskip

  {\bf \small \noindent Step 1.} We first note that
   \eqref{eq:phiz} is the Euler-Lagrange equation for the energy
  functional 
  \begin{align} 
    \label{eq:E}
    {\cal E}[\phi]=\int_\mathbb{R} \Biggl\{ \frac{1}{2} \left({d \phi
        \over d \zeta} \right)^2 + {\eta \over p - 1}
    -\frac{\phi^2(p+1-2\phi^{p-1})}{(p-1)^2} \Biggr\} d \mu,
\end{align}
where $\eta(\zeta)$ is as in the statement of the theorem. Indeed, the
functional $\mathcal E$ in \eqref{eq:E} is continuously differentiable
in $H^1(\mathbb R, d \mu)$ in the natural admissible class $\mathcal
A$ defined as:
\begin{eqnarray}
  \label{eq:A}
  \mathcal A := \{ \phi \in H^1_\mathrm{loc}(\mathbb R): \phi - \eta
  \in H^1(\mathbb R, d \mu), \hspace{1mm} 0 \leq \phi \leq 1 \}.  
\end{eqnarray}
Note that the role of $\eta$ in the definition of $\mathcal E$ is to
ensure that the integral in \eqref{eq:E} converges for all $\phi \in
\mathcal A$. The precise form of $\eta(\zeta)$ is unimportant. Then it
is easy to see that the weak form of \eqref{eq:phiz} in $H^1(\mathbb
R, d \mu)$ is precisely the condition that the Fr\'echet derivative of
$\mathcal E[\phi]$ is zero.

\smallskip

{\bf \small \noindent Step 2.} We now establish weak sequential
lower-semicontinuity and coercivity of the functional $\mathcal E$ in
the admissible class $\mathcal A$ in the following sense: let $\phi_k
= \eta + w_k$, where $w_k \rightharpoonup w$ in $H^1(\mathbb R, d
\mu)$. Then 1) $\liminf_{k \to \infty} \mathcal E[\phi_k] \geq
\mathcal E[\phi]$, where $\phi = \eta + w$, and 2) if $\mathcal
E[\phi_k] \leq M$ for some $M \in \mathbb R$, then
$||w_k||_{H^1(\mathbb R, d \mu)} \leq M'$ for some $M' > 0$.

Let us introduce the notation $\mathcal E[\phi, (a, b)]$ for the
integral in \eqref{eq:E}, in which integration is over all $\zeta \in
(a, b)$.  Then, using \eqref{eq:poinc} from Lemma \ref{l:poinc} we
find that for $R \geq 1$
\begin{eqnarray}
  \label{eq:ERp}
  \mathcal E[\phi_k, (R, +\infty)] \geq \left( e^{2 R} - {p + 1 \over
      (p - 1)^2} \right) \int_R^\infty w_k^2 d \mu > 0.
\end{eqnarray}
Similarly, taking into account that the integrand in \eqref{eq:E} is
non-negative for $\zeta \leq 0$, we have $\mathcal E[\phi_k, (-\infty,
-R)] \geq 0$ for every $R \geq 0$. Since $\mathcal E[\cdot, (-R, R)]$
is lower-semicontinuous by standard theory \cite{dalmaso}, we obtain
$\mathcal E[\phi_k] \geq \mathcal E[\phi_k, (-R, R)]$, yielding the
first claim by passing to the limit $R \to \infty$.

To prove coercivity, we first note that by  \eqref{eq:poinc}
\begin{eqnarray}
  \label{eq:Ec}
  \mathcal E[\phi_k, (R, +\infty)]\geq \int_R^\infty \left\{ \frac12
    \left( { d w_k \over d \zeta} \right)^2 - {p + 1 \over (p
      - 1)^2} w_k^2 \right\}  d \mu \nonumber \\
  \geq \frac14 \int_R^\infty \left\{ \left( { d w_k \over d \zeta}
    \right)^2 + w_k^2 \right\} d \mu,
\end{eqnarray}
for large positive $R$. On the other hand, since $p - 1 - \phi^2 (p +
1 - 2 \phi^{p-1}) \geq (p - 1) (1 - \phi)^2$ for all $0 \leq \phi \leq
1$, we have
\begin{eqnarray}
  \label{eq:Emc}
  \mathcal E[\phi_k, (-\infty, 0)] \geq
  \int_{-\infty}^0 \left\{ \frac12 \left( { d w_k \over d \zeta}
    \right)^2 + {w_k^2 \over p - 1} \right\} d \mu. 
\end{eqnarray}
Finally, by boundedness of $\phi_k$ and $\eta$, we also have
\begin{eqnarray}
  \label{eq:E0c}
  \mathcal
  E[\phi_k, (0, R)] \geq \frac12 \int_0^R \left\{ \left(
      { d w_k \over d \zeta} \right)^2 + w_k^2 \right\} d
  \mu - C R, 
\end{eqnarray}
for some $C > 0$ independent of $w_k$. So the second claim follows.

\smallskip

{\bf \small \noindent Step 3.} In view of the lower-semicontinuity and
coercivity of $\mathcal E$ proved in Step 2, by the direct method of
calculus of variations there exists a minimizer $\phi \in \mathcal A$
of $\mathcal E$. Noting that since the barriers $\phi = 0$ and $\phi =
1$ solve \eqref{eq:phiz} as well, we also have (see
e.g. \cite{kinderlehrer}) that $\phi$ is a weak solution of
\eqref{eq:phiz} by continuous differentiability of $\mathcal E$ in
$H^1(\mathbb R, d \mu)$ noted in Step 1. Furthermore, by standard
elliptic regularity theory \cite{gilbarg}, $\phi \in C^\infty(\mathbb
R)$ and is, in fact, a classical solution of \eqref{eq:phiz}. Also, by
strong maximum principle \cite{gilbarg}, we have $0 < \phi < 1$. To
show monotonicity, suppose, to the contrary, that $\phi(a) < \phi(b)$
for some $a < b$. Then $\phi(\zeta)$ attains a local minimum for some
$\zeta_0 \in (-\infty, b)$. However, by \eqref{eq:phiz} we have $d^2
\phi(\zeta_0) / d \zeta^2 < 0$, giving a contradiction. By the same
argument $d \phi / d \zeta = 0$ is also impossible for any $\zeta \in
\mathbb R$. Finally, since $\phi - \eta \in H^1(\mathbb R, d \mu)$,
monotonicity implies the first condition in \eqref{eq:phibcz1} and
\eqref{eq:phibcz2}.

\smallskip

{\bf \small \noindent Step 4.} We now discuss the asymptotic behavior
of minimizers obtained in Step 3 as $\zeta \to \pm \infty$ and, in
particular, prove the second parts of \eqref{eq:phibcz1} and
\eqref{eq:phibcz2} and the fact that every solution of \eqref{eq:phiz}
belonging to $\mathcal A$ has the same asymptotic decay, which will be
needed later.  Let us first consider the case of $\zeta\to +\infty$.
Performing the Liouville transformation by introducing
 \begin{eqnarray}\label{eq: Lu+}
 \psi
= \phi \sqrt{\rho} \in L^2(R, +\infty), 
\end{eqnarray}
where $\rho$ is defined in
 \eqref{eq:muu} and $R \geq 1$ is arbitrary, we rewrite
 \eqref{eq:phiz} in the form
\begin{eqnarray}
  \label{eq:schr}
  {d^2 \psi \over d \zeta^2} = q(\zeta) \psi, \qquad \zeta \geq R.  
\end{eqnarray}
Here $q(\zeta) = q_0(\zeta) +
q_1(\zeta)$, where
\begin{align}
  \label{eq:q0}
  q_0(\zeta) & = \frac14\left(
    \frac{e^{4\zeta}}{4}+\frac{p-5}{p-1}e^{2\zeta} +1
  \right), \\ 
  \label{eq:q1}
  q_1(\zeta) & = \frac{2(p+1)}{(p-1)^2} \, \phi^{p-1} (\zeta). 
\end{align}
Observe that $q(\zeta) \geq q_0(\zeta) \geq \tfrac14 > 0$ for all
$\zeta \geq R$, with $R$ sufficiently large positive. Therefore,
\eqref{eq:schr} has two linearly-independent positive solutions
$\psi_1$ and $\psi_2$, such that $\psi_1 \to 0$ and $\psi_2 \to
\infty$ together with their derivatives as $\zeta \to +\infty$ (see
e.g. \cite{sansone}). In particular, $\psi = C \psi_1 \in L^2(R,
+\infty)$ for some $0<C<\infty $, and by direct computation
\begin{align}
  \label{eq:phiz0}
  {d \phi \over d \zeta} = {C \over \sqrt{\rho}} \left( {d \psi_1
      \over d \zeta} - \frac{\psi_1}{2} {d \over d \zeta} \ln \rho
  \right) \to 0 \qquad \text{as} \quad \zeta \to +\infty.
\end{align}

On the other hand, as follows from \eqref{eq:poincR}, we have
\begin{eqnarray}
  \label{eq:q11}
  q_1(\zeta) = o( \rho^{1-p \over 2}),
\end{eqnarray}
so $q_1(\zeta)$ has a super-exponential decay as $\zeta \to
+\infty$. Let $\psi_0$ be the unique positive solution of
\eqref{eq:schr} with $q = q_0$ and $\psi_0(R) = 1$ which goes to zero
as $\zeta \to +\infty$. Then we claim that $\psi_1(\zeta) /
\psi_0(\zeta) \to c$ for some $0 < c < \infty$.  Indeed, functions
$\psi_1$ and $\psi_0$ satisfy
 \begin{eqnarray}
  \label{eq:schr01}
  {d^2 \psi_1 \over d \zeta^2} = (q_0(\zeta)+q_1(\zeta)) \psi_1,
  \quad  {d^2 \psi_0 \over d \zeta^2} = q_0(\zeta) \psi_0,\qquad \zeta
  \geq R.   
\end{eqnarray}
 %after
%straightforward algebraic manipulations we have
Multiplying the first and the second equation of \eqref{eq:schr01} by
$\psi_0$ and $\psi_1$, respectively, and taking the difference, we
obtain
\begin{eqnarray}
  \frac{d}{d\zeta}\left( \psi_0 {d \psi_1 \over d \zeta} -  \psi_1{d
      \psi_0 \over d \zeta} \right)=   q_1(\zeta)\psi_0\psi_1.
\end{eqnarray}
Integrating this equation and taking into account that $\psi_0,\psi_1$
and their derivatives vanish as $\zeta\to +\infty$, we have
\begin{eqnarray}
  \psi_0 (\zeta){d \psi_1(\zeta) \over d \zeta} -  \psi_1(\zeta){d
    \psi_0 (\zeta)\over d \zeta} = -\int_{\zeta}^{\infty}
  q_1(s)\psi_0(s)\psi_1(s)ds ,
\end{eqnarray}
and therefore
%\begin{eqnarray}
%\frac{d}{d\zeta}\ln \left( {\psi_1(\zeta) \over \psi_0(\zeta)}\right)=
%-\int_{\zeta}^{\infty} q_1(s)\frac{\psi_0(s)\psi_1(s)}{\psi_0(\zeta)\psi_1(\zeta)}ds
%\end{eqnarray}
\begin{eqnarray}
  \label{eq:log}
  {d \over d \zeta} \ln \left( {\psi_1 \over \psi_0} \right) = -
  \int_\zeta^{\infty} q_1(s) {\psi_1(s) \psi_0(s) \over \psi_1(\zeta)
    \psi_0(\zeta)} ds.  
\end{eqnarray}
Integrating this equation again, we obtain
\begin{eqnarray}
  \label{eq:log1}
  \ln \left( {\psi_1(\zeta) \over \psi_0(\zeta) } \right) =
  \ln \left( {\psi_1(R) \over \psi_0(R) } \right)-
  \int_R^{\zeta}  \int_\sigma^{\infty} q_1(s) {\psi_1(s) \psi_0(s)
    \over \psi_1(\sigma) 
    \psi_0(\sigma)} dsd\sigma.  
\end{eqnarray}
In a view of boundedness of functions $\psi_0$ and $\psi_1$, we have
$|\psi_0(s)/\psi_0(\sigma)|, |\psi_1(s)/\psi_1(\sigma)|\le C$ for some
$C > 0$ and all $s\ge \sigma\ge R$.  Moreover, the estimate in
\eqref{eq:q11} gives $|q_1(s)|\le C' \exp(-s)$ for some $C' > 0$ and
all $s\in[R,\infty)$.  Therefore, the integral in the right-hand side
of \eqref{eq:log1} converges:
 \begin{eqnarray}
   \int_R^{\zeta}  \int_\sigma^{\infty} \left| q_1(s) {\psi_1(s)
       \psi_0(s) \over \psi_1(\sigma) 
       \psi_0(\sigma)} \right| dsd\sigma\le C
   \int_R^{\zeta}  \int_\sigma^{\infty} e^{-s} dsd\sigma  \leq C
   e^{-R}<\infty, 
\end{eqnarray}
which immediately implies that the ratio of $\psi_0$ and $\psi_1$
approaches a finite non-zero limit as $\zeta \to +\infty$.

We can use a similar treatment to establish the asymptotic behavior of
minimizers when $\zeta \to -\infty$.  The Liouville transformation
 \begin{eqnarray}\label{eq: Lu-}
 \theta
= (1-\phi) \sqrt{\rho} \in L^2(-\infty, R), 
\end{eqnarray}
with $\rho$ defined by \eqref{eq:muu} and arbitrary $R \leq 0$ applied
to \eqref{eq:phiz} yields
\begin{eqnarray}
  \label{eq:schr2}
  {d^2 \theta \over d \zeta^2} = r(\zeta) \theta, \qquad \zeta \leq R.  
\end{eqnarray}
Here $r(\zeta) = r_0(\zeta) +
r_1(\zeta)$, where
\begin{align}
  \label{eq:p0}
  r_0(\zeta) & = \frac14\left(\left(\frac{3p+1}{p-1}\right)^2
   +\frac{p-5}{p-1}e^{2\zeta} + \frac{e^{4\zeta}}{4}
  \right), \\ 
  \label{eq:p1}
  r_1(\zeta) & = \frac{2(p+1)}{(p-1)^2} \left(\frac{\phi(1-
      \phi^{p-1})}{1-\phi}+1-p\right) .
\end{align}
By direct computation, note that in the limit $\zeta \to -\infty$ we
have
\begin{align}
  \label{eq:r12}
  r_0(\zeta) \to \frac14\left(\frac{3p+1}{p-1}\right)^2, \qquad
  r_1(\zeta) \to 0^-.
\end{align}
Therefore, $r_0(\zeta) \geq r(\zeta) \geq \tfrac14 > 0$ for all $\zeta
\leq R$ with $R$ sufficiently large negative, and \eqref{eq:schr2} has
two linearly-independent positive solutions $\theta_1$ and $\theta_2$
such that $\theta_1 \to 0$ and $\theta_2 \to \infty$ together with
their derivatives as $\zeta \to -\infty$. In particular, $\theta= C
\theta_1 \in L^2( -\infty,R)$ for some $0<C<\infty$, and 
\begin{align}
  \label{eq:phiz1}
  {d \phi \over d \zeta} = -{C \over \sqrt{\rho}} \left( {d \theta_1
      \over d \zeta} - \frac{\theta_1}{2} {d \over d \zeta} \ln \rho
  \right) \to 0 \qquad \text{as} \quad \zeta \to -\infty.
\end{align}

On the other hand, as follows from \eqref{eq:poincR+} we have
\begin{eqnarray}
  \label{eq:r11}
  r_1(\zeta) = o( \rho^{-1/2}),
\end{eqnarray}
so $r_1(\zeta)$ has an exponential decay as $\zeta \to -\infty$.
Computations practically identical to those presented above show that
the ratio of $\theta_0$ (the solution of \eqref {eq:schr2} with
$r=r_0$ which decays as $\zeta \to-\infty$) and $\theta_1$ tends to a
positive constant as $\zeta \to -\infty$.

%The result then follows by the decay of $\psi_0$ and $\psi_1$ and the
%estimate in  \eqref{eq:q11} upon integration of  \eqref{eq:log}.

\smallskip

{\bf \small \noindent Step 5.} We now prove uniqueness of the obtained
solution, taking advantage of a sort of convexity of $\mathcal E$
similar to the one pointed out in \cite{kawohl85}. Suppose, to the
contrary, that there are two functions $\phi_1, \phi_2 \in \mathcal A$
which solve \eqref{eq:phiz}. Define
\begin{align}
  \label{eq:phit}
  \phi^t := \sqrt{t \phi_2^2 + (1 - t) \phi_1^2}.
\end{align}
We claim that $\phi^t \in \mathcal A$ as
well. Indeed, in view of the result of Step 4 we have $m < \phi_1 /
\phi_2 < M$ for some $M > m > 0$ and, therefore,
\begin{align}
  ||\phi^t||_{L^2((0, 1), d \mu)} & \leq C, \quad
  ||\phi^t||_{L^2((1, \infty), d \mu)}^2 \leq ||\phi_1||_{L^2((1,
    \infty), d \mu)}^2 + ||\phi_2||_{L^2((1, \infty), d \mu)}^2, \\
  ||1 - \phi^t||_{L^2((-\infty, 0), d \mu)}^2 & = \int_{-\infty}^0
  \left( { 1 - t \phi_2^2 - (1 - t) \phi_1^2 \over 1 + \sqrt{t
        \phi_2^2
        + (1 - t) \phi_1^2} } \right)^2 d \mu \nonumber \\
  & \leq C ( ||1 - \phi_1||_{L^2((-\infty, 0), d \mu)} + ||1 -
  \phi_2||_{L^2((-\infty, 0), d \mu)} )^2, \\
  ||d\phi^t/d \zeta||_{L^2(\mathbb R, d \mu)}^2 & = \int_\mathbb{R} {1
    \over t \phi_2^2 + (1 - t) \phi_1^2 } \left( t \phi_2 {d \phi_2
      \over d \zeta} + (1 - t)
    \phi_1 {d \phi_1 \over d \zeta} \right)^2 d \mu \nonumber \\
  & \leq C (|| d \phi_1 / d \zeta||_{L^2(\mathbb R, d \mu)} + ||d
  \phi_2 / d \zeta||_{L^2(\mathbb R, d \mu)})^2,
\end{align}
for some $C > 0$.  In fact, it is easy to see that the function $E(t)
:= \mathcal E[\phi^t]$ is twice continuously differentiable for all $t
\in [0, 1]$. A direct computation yields
\begin{eqnarray}
  \label{eq:dE2}
  {d^2 E(t) \over d t^2} = \int_\mathbb{R} \Biggl\{ {\phi_1^2 \phi_2^2
    \over (t \phi_2^2 + (1 - t) \phi_1^2)^3 } \left( \phi_2 {d \phi_1
      \over d \zeta} - \phi_1 {d \phi_2 \over d \zeta}\right)^2
  \nonumber \\
  + {p+ 1 \over 2 p - 2} (\phi_1^2 - \phi_2^2)^2  (t \phi_2^2 + (1
  - t) \phi_1^2)^{p - 3 \over 2} \Biggr\} d \mu(\zeta). 
\end{eqnarray}
Therefore, $d^2 E(t) / dt^2 > 0$ for all $t \in [0,1]$, and so $E(t)$
is strictly convex. However, since the map $t \mapsto \phi^t - \eta$
is of class $C^1([0, 1]; H^1(\mathbb R, d \mu))$, which can be seen by
a computation analogous to the one in \eqref{eq:dE2}, this contradicts
the fact that $dE(0)/dt = dE(1)/dt = 0$ by the assumption that
$\phi_1$ and $\phi_2$ solve weakly \eqref{eq:phiz} and hence are
critical points of $\mathcal E$. \hfill $\Box$

\begin{remark}
  Results of Step 4 of the proof above allow to obtain the precise
  asymptotic behavior of the solution of \eqref{eq:phiz} constructed
  in Theorem \ref{t:exun} by using the exact solutions of the
  associated linearizations of \eqref{eq:phiz} about $\phi=0$ and
  $\phi=1$.  These asymptotics read \cite{mgs:pre11}:
\begin{align} \label{eq:rem1}
  \phi(\zeta)\sim\exp\left(-\frac{e^{2\zeta}}{4}+\frac{5-p}{p-1}\zeta\right),
  \quad & \zeta\to +\infty,\nonumber \\ 
  1-\phi(\zeta)\sim \exp\left(\frac{2 (p+1)}{p-1}\zeta\right),\quad
  & \zeta \to -\infty. 
\end{align}

\end{remark}

\section{Long time behavior of solutions for problem \eqref{eq:u}} 
\label{s:long}

In this section we prove that the ultra-singular solutions constructed
in Sec.~\ref{sec:similarity-ansatz} have a direct relevance to the
long time behavior of solutions for the problem in
\eqref{eq:u}. Specifically, solutions of \eqref{eq:u} converge to
self-similar profile $\phi$ at the fixed ratio $x/\sqrt{t}$ as $t\to
\infty$. That is, the following result holds:
\begin{theorem}
  \label{t:conv}
  Given $\alpha>0$, let $u$ and $v$ be the solutions of \eqref{eq:u}
  and \eqref{eq:v}, respectively, and set
 \begin{eqnarray}
F(\zeta,t)=\frac{u(x,t)}{v(x)}, \quad \zeta=\ln\left(\frac{x}{\sqrt{t}}\right).
\end{eqnarray}
Then
 \begin{eqnarray}
   \lim_{t\to\infty} F(\zeta,t)=\phi(\zeta) \quad \forall \zeta\in
   \mathbb{R}. 
 \end{eqnarray}
  Moreover,
  \begin{eqnarray}
  \phi(\xi) \le F(\zeta,t)\le \phi(\zeta)
  \end{eqnarray}
  where $\xi(\zeta,t)=\ln(e^{\zeta}+b t^{-1/2})$ and $b$ is some large
  enough constant.
  \end{theorem}

  {\em \bf \small \noindent Proof. } The proof relies on a direct
  application of the comparison principle.  We start with a
  formulation of the comparison principle which will be applied to
  \eqref{eq:u}.  Define the following quantities
\begin{eqnarray}
&& P[u]=u_t-u_{xx}+u^p,\\
&& Q[u]=u_x+\alpha,
\end{eqnarray}
assume that the functions $\bar u$ and $\underline u$ satisfy the
differential inequalities
%For super-solution $\bar u$ we have
\begin{eqnarray}\label{eq:super1}
  && P[\bar u]\ge 0, \quad t>0, \quad x>0, \\
  && Q[\bar u]\le 0, \quad t>0, \quad x=0, \label{eq:super2}
\end{eqnarray}
and
%For sub-solution $\underline u$ we have
\begin{eqnarray}\label{eq:sub1}
  && P[\underline u]\le 0, \quad t>0, \quad x>0, \\
  && Q[\underline u]\ge 0, \quad t>0, \quad x=0. \label{eq:sub2}
\end{eqnarray}
and, in addition, assume that $\bar u(x,t=0)=\underline u(x,t=0)=0$.
Such functions are called super- and sub-solutions for \eqref{eq:u}
and have the property \cite{protter}:
\begin{eqnarray}
  \label{eq:comp}
  \underline u(x,t) \le u(x,t) \le \bar u(x,t), \quad
  (x,t)\in[0,\infty)\times[0,\infty).
\end{eqnarray}
In what follows we will explicitly construct sub- and super-solutions
for \eqref{eq:u}.

We first show that the function
\begin{eqnarray}\label{eq: usub}
\underline u(x,t)=v(x)\phi(z), \quad z=\ln\left(\frac{x+b}{\sqrt{t}}\right),
\end{eqnarray}
is a sub-solution, provided $b\ge a$ is large enough. Here $\phi$
verifies \eqref{eq:phiz}, \eqref{eq:phibcz1} and \eqref{eq:phibcz2},
and $a$ is defined in \eqref{eq:a}.
% \begin{eqnarray} \label{eq: phi}
%   \phi^{\prime\prime}+\left(\frac{e^{2z}}{2}-\frac{n+3}{n-1}\right)\phi^{\prime}+2\frac{(n+1)}{(n-1)^2}\phi(1-\phi^{n-1})=0, 
%   \qquad z\in(-\infty,\infty).
%\end{eqnarray}
Direct substitution of \eqref{eq: usub} into \eqref{eq:sub1} gives:
\begin{align}
  P[\underline u] & =\frac{v(x)}{(x+b)^2} \nonumber \\ & \times
  \left\{ \frac{4}{p-1}\left(1-\frac{x+b}{x+a}\right)
    \left(-\frac{d}{dz}\phi\right)
    +\frac{2(p+1)}{(p-1)^2}\left(1-\left(\frac{x+b}{x+a}\right)^2\right)
    \phi(1-\phi^{p-1})\right\}.
\end{align}
In view of the fact that $d\phi/dz < 0$ we have
\begin{eqnarray}
  P[\underline u]\le 0 \qquad \forall x>0, \quad \forall t>0,
\end{eqnarray}
provided that $b\ge a$.  

Next, direct computations also give
\begin{eqnarray}
  Q[\underline u(x=0,t)]=\frac{2A}{(p-1)a^{\frac{p+1}{p-1}}}\left(
    1-\phi(z_b)     +\frac{1}{b}\frac{a(p-1)}{2}
    \frac{d}{dz}\phi(z_b)\right), \quad
  z_b=\ln\left(\frac{b}{\sqrt{t}}\right) .
\end{eqnarray}
Let us show that $Q[\underline u(x=0,t)]\ge 0$ for $t>0$ when $b$ is
large.  To do so, it is enough to show that
\begin{eqnarray}
  \label{gzp}
  g(z) :=1-\phi(z)+\eps \frac{d}{dz}\phi(z)\ge 0 \qquad \forall z\in
  \mathbb{R}, 
\end{eqnarray}
for $\eps>0$ small.  Indeed, observe first that $\lim_{z \to +\infty}
g(z) = 1$ and $\lim_{z \to -\infty} g(z) = 0$.  So, if \eqref{gzp} is
violated, $g(z)$ has a local minimum at some point $z^* \in \mathbb R$
with $g(z^*)<0$. Since $z^*$ is a critical point we have
\begin{eqnarray}
  0 &=
  &\frac{d}{dz}g(z^*)=-\frac{d}{dz}\phi(z^*)+\eps\frac{d^2}{dz^2}\phi(z^*)=
  \nonumber \\ 
  &&-\left(1+\eps
    \frac{e^{2z^*}}{2}-\eps \,
    \frac{p+3}{p-1}\right)\frac{d}{dz}\phi(z^*)-2\eps 
  \frac{(p+1)}{(p-1)^2}(\phi(z^*)-\phi^{p}(z^*)). 
\end{eqnarray}
Therefore, there exists $\eps \in (0, 1)$ such that
\begin{eqnarray}
  |\phi_z(z^*)|\le 1-\phi(z^*),
\end{eqnarray}
Thus, from the definition of $g$ we have
\begin{eqnarray}
  g(z^*)\ge \left(1-\eps \right) (1-\phi(z^*)) \geq 0,
\end{eqnarray}
contradicting our assumption about $g(z^*)$.  Finally, choosing
$b=\max\{a,\frac{a(p-1)}{2\eps}\}$ we have that the conditions in
\eqref{eq:sub1} and \eqref{eq:sub2} are satisfied and thus \eqref{eq:
  usub} is indeed a sub-solution for $u$.

Now we turn to the construction of a super-solution, which we will
seek in the form
\begin{eqnarray}\label{eq: usup}
  \bar u(x,t)=v(x)\phi(\zeta), \quad
  \zeta=\ln\left(\frac{x}{\sqrt{t}}\right). 
\end{eqnarray}
Straightforward computations give
\begin{align}
  P[\bar u(x,t)]= & \frac{v(x)}{x^2}\Big\{
  \frac{4}{p-1}\left(1-\frac{x}{x+a}\right)\left(-\frac{d}{d\zeta}\phi\right)
  \nonumber \\
  & + \frac{2(p+1)}{(p-1)^2}\left(1-\left(\frac{x}{x+a}\right)^2
  \right)\phi(1-\phi^{p-1})\Big\} ,
\end{align}
and
\begin{eqnarray}
  Q[\bar u(x=0, t)]=\frac{2A}{(p-1)a^{\frac{p+1}{p-1}}}
  \lim_{\zeta \to-\infty} \left( 1-\phi(\zeta)
    +\frac{1}{\sqrt{t}}\frac{a(p-1)}{2}e^{-\zeta}\frac{d}{d\zeta}\phi(\zeta)\right),
\end{eqnarray}
It is clear that $P[\bar u(x,t)]\ge 0$ for all $t>0$ and $x>0$.  Let
us now show that
\begin{eqnarray}
  Q[\bar u(x=0,t)]=0 \qquad \forall t>0. 
\end{eqnarray}
Since by \eqref{eq:phibcz1} and \eqref{eq:phibcz2}
\begin{eqnarray}
\lim_{\zeta\to-\infty}(1-\phi(\zeta))=0,
\end{eqnarray}
we only need to show that
\begin{eqnarray}\label{eq:epz}
 \lim_{\zeta\to-\infty} e^{-\zeta}\frac{d}{d\zeta} \phi(\zeta) =0.
 \end{eqnarray}
Indeed, multiplying
\eqref{eq:phiz} by $\rho$ we have
\begin{eqnarray}
  \frac{d}{d\zeta}\left( \rho
    \frac{d}{d\zeta}\phi\right)=-\frac{2(p+1)}{(p-1)^2}\rho\phi(1-\phi^{p-1}).  
\end{eqnarray}
Integrating this equation and rearranging terms involving $\rho$, we
obtain
\begin{eqnarray}
  \label{eq:limit+}
  e^{-\zeta}\frac{d}{d\zeta}\phi(\zeta) & = & \exp\left(-\frac{e^{2\zeta}}{4}
    +\frac{4}{p-1}\zeta\right) 
  \nonumber \\ 
  && \times \left(\rho(R)\frac{d}{d\zeta}\phi(R)
    +\frac{2(p+1)}{(p-1)^2}\int_\zeta^R\rho(s)\phi(s)
    (1-\phi^{p-1}(s))ds\right). 
\end{eqnarray}
By \eqref{eq:rem1} we have $\rho(\zeta)\phi(\zeta)
(1-\phi^{p-1}(\zeta))\sim \exp(\zeta)$ as $\zeta\to-\infty$ and thus
the integral in the right-hand side of \eqref{eq:limit+} converges as
$\zeta\to-\infty$, which readily implies \eqref{eq:epz}.  Therefore,
both conditions \eqref{eq:super1} and \eqref{eq:super2} are satisfied
and so \eqref{eq: usup} is a super-solution.

Finally, the statement of the theorem follows from \eqref{eq:comp},
\eqref{eq: usub} and \eqref{eq: usup}. \hfill $\Box$

\begin{remark}
 \label{r:2}
 Note that the result of Theorem \ref{t:conv} may be extended to
 problem \eqref{eq:u} in which the constant $\alpha$ is replaced by a
 bounded, monotonically increasing function $\alpha(t) > 0$.
\end{remark}

\paragraph{Acknowledgements.} This work was supported, in part, by NSF
via grant DMS-1119724. CBM would also like to acknowledge partial
support by NSF via grant DMS-0908279. We wish to thank S. Shvartsman
for suggesting this problem to us and V. Moroz for helpful
comments. PVG also would like to acknowledge valuable discussions with
S. Kamin.

\bibliographystyle{plain}

\bibliography{../mura,../nonlin,../bio,../egfr}

\end{document}